\documentclass{article}
\usepackage{epsf}
\usepackage{fullpage}
\usepackage{euscript,amsmath,amssymb, amsthm} %, amscd}
\usepackage[dvips]{graphicx}
\usepackage[dvips]{color}
\renewcommand{\k}{\mathbf{k}}

\newcommand{\Z}{\mathbf{Z}}

\newcommand{\dm}{\operatorname{\mathsf{dim}}}

\renewcommand{\span}{\operatorname{\mathsf{span}}}

\newcommand{\Id}{\operatorname{\mathsf{Id}}}

\newcommand{\iA}{\mathsf{A}}
\newcommand{\iB}{\mathsf{B}}

\newcommand{\s}{\mathsf{s}}

\newcommand{\X}{\mathsf{X}}
\newcommand{\y}{\mathsf{y}}

\newcommand{\impl}{\Longrightarrow}

\newcommand{\nozero}{\setminus \{0\}}

\newcommand{\defeq}{\stackrel{\mathrm{def}}{=}}

\def\cent#1{\mathcal{Z}(#1)}

\newcommand{\Hh}{\mathbf{H}}

\newcommand{\subH}{\overline{\Hh}}
\newcommand{\subsubH}{\overline{\overline{\Hh}}}

\newcommand{\pf}{\mathbf{F}_p}
\newcommand{\vx}{V_{\X}}
\newcommand{\vgen}{V_{\text{gen}}} 

\begin{document}

\newtheorem{theorem}{Theorem}[section]
\newtheorem{proposition}[theorem]{Proposition}
\newtheorem{lemma}[theorem]{Lemma}
\newtheorem{claim}[theorem]{Claim}  
\newtheorem{corollary}[theorem]{Corollary}

\theoremstyle{definition}
\newtheorem{definition}[theorem]{Definition}
\newtheorem{definition-sort-of}[theorem]{`Definition'}
\newtheorem{notation}[theorem]{Notation}

\newtheorem{example}[theorem]{Example}
\newtheorem{remark}[theorem]{Remark}
\newtheorem{fact}[theorem]{Fact} 
\newenvironment{block}[1]{\medskip \noindent {\Large \bf #1.}}{\medskip}

\title{Representations of trigonometric Cherednik algebras of rank 1
in positive characteristic}

\author{Fr\'ed\'eric Latour}

\maketitle

\section{Introduction}

Cherednik's double affine Hecke algebras are an important class of
algebras attached to root systems. They were introduced in \cite{Ch3} as a 
tool of
proving Macdonald's conjectures, but are also interesting by themselves, since 
they provide universal deformations of twisted group algebras of double affine 
Weyl groups. One may distinguish rational, trigonometric, and elliptic Cherednik
algebras, which contain 0, 1, and 2 copies of the root lattice, respectively 
(rational and trigonometric algebras are degenerations of the elliptic ones; see 
\cite{EG}).

Development of representation theory of Cherednik algebras (in particular, 
description of all irreducible finite dimensional representations) is an 
important open problem. In the characteristic zero case, it is solved completely 
only for type A, while in other types only partial results are available (see
\cite{EG},\cite{BEG}, and \cite{ChO} for the rank 1 case). In positive
characteristic, the rank 1 case (in the more general setting of complex 
reflection groups) is settled by the author in \cite{L}, after which the higher 
rank case (of type A) was considered in \cite{FG}.

The goal of this paper is to extend the results of \cite{L} to
the trigonometric case. That is, we study the representation theory
of trigonometric Cherednik algebras in positive characteristic $p$
in the simplest case of rank 1.
Our main result is a complete description of irreducible representations
of such algebras.

The paper is organized as follows.

In Section 2, we state the main results.

In Section 3, we prove the results for the ``classical'' case, i.e.
the case when the ``Planck's constant'' $t$ is zero. In this case,
generic irreducible representations have dimension $2$; one-dimensional
representations exist when the ``coupling constant'' $k$ is zero.

In Section 4, we prove the results for the ``quantum'' case, i.e.
the case when the ``Planck's constant'' $t$ is nonzero. In this case,
generic irreducible representations have dimension $2p$; smaller
representations exist when the ``coupling constant'' $k$ is an element of
$\pf \subset \k$; namely, if $k$ is an integer with $0 \leq k \leq p-1$,
then there exist irreducible representations of dimensions $p-k$ and
$p+k$.

{\bf Acknowledgements.} The author thanks his adviser Pavel Etingof for
posing the problem and useful discussions, as well as for helping to write
an introduction. The work of the author
was partially supported by the National Science Foundation (NSF) grant
DMS-9988796 and by a Natural Sciences and Engineering Research
Council of Canada (NSERC) Julie Payette research scholarship.

\section{Statement of Results}

Let $\k$ be an algebraically closed field of characteristic $p$, where $p 
\neq 2$. Let $t, k \in \k$, and
let $\Hh(t,k)$ be the algebra (over $\k$) generated by $\X, 
\X^{-1}, 
\s$ and $\y$, subject to the following relations:

\begin{eqnarray}
\s \X &=& \X^{-1} \s \label{rel1} \\
\s^2 &=& 1 \label{rel2} \\
\s \y + \y \s &=& -k \label{rel3} \\
\X \y \X^{-1} &=& \y - t + k\s. \label{rel4}
\end{eqnarray}

We will classify the irreducible representations of $\Hh(t,k)$.
Now, for $t \neq 0$, $\Hh(t,k)$ is clearly isomorphic to 
$\Hh(1,\frac{k}{t})$ under the 
map $$\X \mapsto \X, \s \mapsto \s, \y 
\mapsto \frac{1}{t} \y.$$ Thus it is sufficient to classify irreducible 
representations of $\Hh(0,k)$ and $\Hh(1,k)$. For brevity we will use the 
notation $\Hh_0 \defeq \Hh(0,k)$ and $\Hh_1 \defeq \Hh(1,k)$, assuming that $k$ 
has been fixed once and for all.

\subsection{Irreducible representations of $\Hh_0$}

\begin{proposition}
Let $k \neq 0$. Then the irreducible representations of $\Hh_0$ are the 
following:
\begin{itemize}
\item For $a, \beta \in \k$, $a, \beta \neq 0$, we have a two-dimensional 
representation $V_{0,1}^{\beta, a}$ with basis $\{ v_0, v_1 \}$, defined by 
the following:
\begin{eqnarray*}
\y v_0 &=& \beta v_0, \\
\y v_1 &=& - \beta v_1, \\
\X v_0 &=& a v_0 - \frac{k^2}{4 \beta^2} v_1, \\
\X v_1 &=& v_0 + 
	\left(\frac{1}{a} - \frac{k^2}{4 a \beta^2} \right) v_1, \\
\s v_0 &=& -\frac{k}{2 \beta} v_0 + 
	\frac{k^3 - 4 k \beta^2}{8a \beta^3} v_1, \\
\s v_1 &=& - \frac{2 a \beta}{k} v_0 + \frac{k}{2 \beta} v_1;
\end{eqnarray*}
\item For $a = \pm 1, b \in \k,$ we have a two-dimensional representation
$V_{0,2}^{a,b}$ with basis $\{ v_0, v_1 \}$, defined by
the following:
\begin{eqnarray*}
\y v_0 &=& 0, \\
\y v_1 &=& v_0, \\
\s v_0 &=& v_0 - k v_1, \\
\s v_1 &=& -v_1, \\
\X v_0 &=& a (v_0 - k v_1), \\
\X v_1 &=& b v_0 + (a - k b) v_1.
\end{eqnarray*}
\end{itemize}
$V_{0,1}^{\beta, a}$ and $V_{0,1}^{\beta', a'}$ are isomorphic if and only 
if
$\beta' = \beta, a' = a$ or
$\beta' = -\beta, a' = \frac{4 \beta^2 - k^2}{4 a \beta^2}.$
$V_{0,2}^{a,b}$ and $V_{0,2}^{a',b'}$ are 
isomorphic if and only if $a = a'$ and $b = b'.$
Furthermore, representations with different subscripts are never 
isomorphic.
\label{prop1}
\end{proposition}

\begin{proposition}
Let $k = 0$. Then the irreducible representations of $\Hh_0$ are the
following:
\begin{itemize}
\item For $a, \beta \in \k$, $a, \beta \neq 0$, we have a two-dimensional 
representation $V_{0,3}^{\beta, a}$ with basis $\{ v_0, v_1 \}$, defined by
the following:
\begin{eqnarray*}
\y v_0 &=& \beta v_0, \\
\y v_1 &=& - \beta v_1, \\
\X v_0 &=& a v_0, \\
\X v_1 &=& \frac{1}{a} v_1, \\
\s v_0 &=& v_1, \\
\s v_1 &=& v_0; 
\end{eqnarray*}
\item For $a \in \k$, $a \notin \{ 0, \pm 1 \},$ we have a two-dimensional 
representation $V_{0,4}^a$ with basis $\{ v_0, v_1 \}$, 
defined by the following:
\begin{eqnarray*}
\y v_0 &=& 0, \\
\y v_1 &=& 0, \\
\X v_0 &=& a v_0, \\
\X v_1 &=& \frac{1}{a} v_1, \\
\s v_0 &=& v_1, \\
\s v_1 &=& v_0;
\end{eqnarray*}
\item For $a = \pm 1, b = \pm 1,$ we have a 
one-dimensional representation $V_{0,5}^{a,b}$
on which $\y, \X$ and $\s$ act as $0, a$ and $b$ 
respectively.
\end{itemize}
$V_{0,3}^{\beta, a}$ and $V_{0,3}^{\beta', a'}$ are isomorphic if and only 
if
$\beta' = \beta, a' = a$ or $\beta' = -\beta, a' = \frac{1}{a}.$
$V_{0,4}^a$ and $V_{0,4}^{a'}$ are isomorphic if and only if $a' = a$ or 
$a' = 
\frac{1}{a}$. $V_{0,5}^{a,b}$ and $V_{0,5}^{a',b'}$ are isomorphic if and 
only if
$a'=a, b'=b.$
Furthermore, representations with different subscripts are never
isomorphic.
\label{prop2}
\end{proposition}

\subsection{Irreducible representations of $\Hh_1$}

\begin{proposition}
Let $k \notin \pf$. Then the irreducible representations of $\Hh_1$ are the
following:
\begin{itemize}
\item For $\mu,d \in \k, d \neq 0, b = (\mu^p - \mu)^2$ with $\frac{k}{2}$ 
not a root of $f(y) = (y^p - y)^2 - b,$ and also for $\mu = 
\pm \frac{k}{2}, d \neq 0,$ 
we have a $2p$-dimensional representation $V_{1,1}^{\mu,d}$
with basis $\{ v_{\mu+j}, v_{-\mu+j}, 
j = 0,1,\ldots, p-1 \}$, defined by the following:
\begin{eqnarray}
\y v_{\beta} &=& \beta v_{\beta}, \quad \beta = \pm \mu, \pm\mu +1, 
\ldots,
\pm \mu + p -1; \label{v11first} \\
\s v_{-\mu-j} &=& - \frac{1}{\mu+j} v_{\mu+j} + \frac{k}{2(\mu+j)}
v_{-\mu-j}, \quad  j = 1, 2, \ldots, p-1; \\
\s v_{\mu+j} &=& \left( \frac{k^2}{4(\mu+j)} - (\mu+j) \right) v_{-\mu-j} - 
\frac{k}{2(\mu+j)} v_{\mu+j}, \quad  j = 1, 2, \ldots, p-1; \\
\s v_{-\mu} &=& \frac{k}{2 \mu} v_{-\mu} - \frac{d}{\mu} v_{\mu}; \\
\s v_{\mu} &=& \left( \frac{k^2}{4d \mu} - \frac{\mu}{d} \right) v_{-\mu} -
\frac{k}{2 \mu} v_{\mu}; \\
\X v_{\beta} &=& \s v_{-\beta-1}, \quad
\beta = \pm \mu, \pm\mu +1, \ldots, \pm \mu + p -1 \label{v11last};
\end{eqnarray}

\item For $\theta = \pm 1,$ we have a $2p$-dimensional representation 
$V_{1,2}^{\theta}$ with basis $\{ v_{j}, w_{j}, j = 0,1,\ldots, p-1 \},$ 
defined by the following:
\begin{eqnarray}
\y v_j &=& j v_j, \quad j = 0, 1, \ldots, p-1; \label{v12first} \\
\y w_j &=& j w_j + v_j, \quad j = 0, 1, \ldots, p-1; \\
\s v_0 &=& -k w_0; \\
\s w_0 &=& -\frac{1}{k} v_0; \\
\s v_{-j} &=& \frac{1}{j} v_j + \frac{k}{2j} v_{-j}, \quad j = 1,
2, \ldots, \frac{p-1}{2}; \\
\s v_{j} &=& \left( j - \frac{k^2}{4j} \right) v_{-j} - \frac{k}{2j} v_{j},
\quad j = 1, 2, \ldots, \frac{p-1}{2}; \\
\s w_{-j} &=& \frac{1}{j^2} v_j + \frac{k}{2j^2} v_{-j} - \frac{1}{j} w_j
+ \frac{k}{2j} w_{-j} \quad j = 1, 2, \ldots, \frac{p-1}{2}; \\
\s w_{j} &=& \left( 1 + \frac{k^2}{4j^2} \right) v_{-j} + \frac{k}{2j^2} v_j -
\frac{k}{2j} w_j + \left( \frac{k^2}{4j} - j \right) w_{-j},
\quad j = 1, 2, \ldots, \frac{p-1}{2}; \\
\X v_{j} &=&  - \s v_{-j-1}, \quad j \neq \frac{p-1}{2}; \\
\X v_{\frac{p-1}{2}} &=& \theta \s
v_{\frac{p-1}{2}}; \label{v12-2ndlast}\\
\X w_{j} &=& \s w_{-j-1}, \quad j \neq \frac{p-1}{2}; \\
\X w_{\frac{p-1}{2}} &=& -\theta \s
w_{\frac{p-1}{2}}. \label{v12last}
\end{eqnarray}
\end{itemize}
$V_{1,1}^{\mu,d}$ and $V_{1,1}^{\mu',d'}$ are isomorphic if
and only if
$$(\mu' - \mu \in \pf \mbox{ {\rm and} } d' = d) \mbox{ {\rm or} }
(\mu' + \mu \in \pf \mbox{ {\rm and} } dd' = \prod_{c \in \pf} 
\left(\frac{k^2}{4} - (\mu+c)^2 \right) ).$$
$V_{1,2}^{\theta}$ and $V_{1,2}^{\theta'}$ are isomorphic if and only if
$\theta = \theta'$
Furthermore, representations with different subscripts are never
isomorphic.
\label{prop3}
\end{proposition}
                                  
Now, in the case where $k \in \pf$, note that there is an isomorphism 
between $\Hh(1,k)$ and $\Hh(1,-k)$, given by
$$\y \mapsto \y, \s \mapsto -\s, \X \mapsto \X, k \mapsto -k.$$
So we may assume that $k$ is an {\em even integer} with $0 \leq 
k \leq p-1.$

\begin{proposition}
Let $k$ be even with $2 \leq k \leq p-1$. Then the irreducible 
representations of $\Hh_1$ are the following:
\begin{itemize}
\item For $\mu,d \in \k, d \neq 0,$ % b = (\mu^p - \mu)^2$, 
we have $V_{1,1}^{\mu,d}$, defined as in Proposition \ref{prop3}. 
\item For $\theta = \pm 1$, we have a 
$(p-k)$-dimensional representation $V_{1,3}^{\theta}$ with basis $\{ 
v_{\frac{k}{2}}, 
v_{\frac{k}{2}+1}, \ldots, v_{-\frac{k}{2}-1} \}$, defined by
\begin{eqnarray}
\y v_j &=& j v_j, \quad j = \frac{k}{2}, \frac{k}{2} + 1 \ldots,
- \frac{k}{2} - 1; \label{v13first} \\
\s v_{-j} &=& \frac{k}{2j} v_{-j} - \frac{1}{j} v_j, \quad
j = \frac{k}{2} + 1, \ldots, \frac{p-1}{2}; \\
\s v_j &=& -j v_{-j} + \frac{k^2}{4j} v_{-j} - \frac{k}{2j} v_j,
\quad j = \frac{k}{2} + 1, \ldots, \frac{p-1}{2}; \\
\s v_{\frac{k}{2}} &=&  -v_{\frac{k}{2}};
\end{eqnarray}
\begin{eqnarray}
\X v_{j} &=& - \s v_{-j-1}, \quad j \neq \frac{p-1}{2}; \\
\X v_{\frac{p-1}{2}} &=& \theta \s
v_{\frac{p-1}{2}} \label{v13last}.
\end{eqnarray}
\item
For $\theta = \pm 1$, we have a $(p+k)$-dimensional representation 
$V_{1,4}^{\theta}$ with basis $\{ v_j, w_i, 
j = 0, \ldots, p-1, i = -\frac{k}{2}, \ldots, \frac{k}{2} - 1 \}$, defined by 
\begin{eqnarray}
\y v_j &=& j v_j, \quad j = 0, 1, \ldots, p-1; \label{v14first} \\
\y w_j &=& j w_j + v_j, \quad j = -\frac{k}{2}, -\frac{k}{2} + 1, \ldots,
\frac{k}{2} - 1; \\
\s v_0 &=& -k w_0; \\
\s w_0 &=& -\frac{1}{k} v_0; \\
\s v_{-j} &=& \frac{1}{j} v_j + \frac{k}{2j} v_{-j},
\quad j = 1, \ldots, \frac{k}{2} - 1, \frac{k}{2} + 1, \frac{p-1}{2}; \\
\s v_{j} &=& \left( j - \frac{k^2}{4j} \right) v_{-j} - \frac{k}{2j} v_{j}, 
\quad j = 1, \ldots, \frac{k}{2} - 1, \frac{k}{2} + 1, \frac{p-1}{2}; \\
\s v_{-\frac{k}{2}} &=& v_{-\frac{k}{2}}; \\
\s v_{\frac{k}{2}} &=& 2 v_{-\frac{k}{2}} - v_{\frac{k}{2}};  \\
\s w_{-j} &=& \frac{1}{j^2} v_j + \frac{k}{2j^2} v_{-j} - \frac{1}{j} w_j
+ \frac{k}{2j} w_{-j} \quad j = 1, \ldots, \frac{k}{2} - 1; \\
\s w_{j} &=& \left( 1 + \frac{k^2}{4j^2} \right) v_{-j}
+ \frac{k}{2j^2} v_j - \frac{k}{2j} w_j + \left( \frac{k^2}{4j} - j \right) 
w_{-j},
\quad j = 1, \ldots, \frac{k}{2} - 1; \\
\s w_{-\frac{k}{2}} &=& -\frac{2}{k} v_{\frac{k}{2}} + \frac{2}{k}
v_{-\frac{k}{2}} + w_{-\frac{k}{2}}; \\
\X v_{j} &=& - \s v_{-j-1}, \quad j \neq \frac{p-1}{2}; \\
\X v_{\frac{p-1}{2}} &=& \theta \s
v_{\frac{p-1}{2}}; \\
\X w_{j} &=& = \s w_{-j-1}, \quad j = -\frac{k}{2}, \ldots,
\frac{k}{2} - 1. \label{v14last}
\end{eqnarray}

\item For $c \in \k$, we have a $2p$-dimensional representation $V_{1,5}^{c}$ 
with basis $\{ v_j, w_i, u_l, j = 0, \ldots, p-1, 
i = -\frac{k}{2}, \ldots, \frac{k}{2} - 1,
l = \frac{k}{2}, \ldots, -\frac{k}{2} -1 \},$ defined by
\begin{eqnarray}
\y v_j &=& j v_j, \quad j = 0, 1, \ldots, p-1; \label{v15first} \\
\y w_j &=& j w_j + v_j, \quad j = -\frac{k}{2}, -\frac{k}{2} + 1, \ldots,
\frac{k}{2} - 1; \\
\s v_0 &=& w_0; \\
\s w_0 &=& v_0; \\
\s v_{-j} &=& \frac{1}{j} v_j + \frac{k}{2j} v_{-j},
\quad j = 1, \ldots, \frac{k}{2} - 1, \frac{k}{2} + 1, \frac{p-1}{2}; \\
\s v_{j} &=& \left( j - \frac{k^2}{4j} \right) v_{-j} - \frac{k}{2j} v_{j},
\quad j = 1, \ldots, \frac{k}{2} - 1, \frac{k}{2} + 1, \frac{p-1}{2}; \\
\s v_{-\frac{k}{2}} &=& v_{-\frac{k}{2}}; \\
\s v_{\frac{k}{2}} &=& 2 v_{-\frac{k}{2}} - v_{\frac{k}{2}};
\end{eqnarray}
\begin{eqnarray}
\s w_{-j} &=& \frac{1}{j^2} v_j + \frac{k}{2j^2} v_{-j} - \frac{1}{j} w_j
+ \frac{k}{2j} w_{-j} \quad j = 1, \ldots, \frac{k}{2} - 1; \\
\s w_{j} &=& \left( 1 + \frac{k^2}{4j^2} \right) v_{-j}
+ \frac{k}{2j^2} v_j  -
\frac{k}{2j} w_j + \left( \frac{k^2}{4j} - j \right) w_{-j},
\quad j = 1, \ldots, \frac{k}{2} - 1; \\
\s w_{-\frac{k}{2}} &=& -\frac{2}{k} v_{\frac{k}{2}} + \frac{2}{k}
v_{-\frac{k}{2}} + w_{-\frac{k}{2}}; \\
\s u_j &=& -\frac{1}{j} u_{-j} - \frac{k}{2j} u_j, \quad j = \frac{k}{2} + 
1,
\ldots, \frac{p-1}{2}; \\
\s u_{-j} &=& \left( \frac{k^2}{4j} -j \right) u_j + \frac{k}{2j} u_{-j}, 
\quad j = \frac{k}{2} + 1, \ldots, \frac{p-1}{2}; \\
\s u_{\frac{k}{2}} &=& \frac{2c}{k} v_{-\frac{k}{2}} - u_{\frac{k}{2}}; \\
\X v_{j} &=& - \s v_{-j-1}, \quad j \neq \frac{p-1}{2}; \\
\X v_{\frac{p-1}{2}} &=& \s u_{\frac{p-1}{2}}; \\
\X u_{j} &=& - \s u_{-j-1}, \quad j = \frac{k}{2}, \ldots,
\frac{p-3}{2}, \frac{p+1}{2}, \ldots, -\frac{k}{2}-1; \\
\X u_{\frac{p-1}{2}} &=& \s v_{\frac{p-1}{2}}.
\label{v15last}
\end{eqnarray}
\end{itemize}
$V_{1,1}^{\mu,d}$ and $V_{1,1}^{\mu',d'}$ are isomorphic if
and only if
$$(\mu' - \mu \in \pf \mbox{ {\rm and} } d' = d) \mbox{ {\rm or} }
(\mu' + \mu \in \pf \mbox{ {\rm and} } dd' = \prod_{c \in \pf} 
\left(\frac{k^2}{4} - (\mu+c)^2 \right)).$$
$V_{1,3}^{\theta}$ and $V_{1,3}^{\theta'}$ are isomorphic if and only if
$\theta = \theta'$.
$V_{1,4}^{c}$ and $V_{1,4}^{c'}$ are isomorphic if and only if
$c = c'$.
$V_{1,5}^{c}$ and $V_{1,5}^{c'}$ are isomorphic if and only if
$c = c'$.
Furthermore, representations with different subscripts are never
isomorphic.
\label{prop4}
\end{proposition}
                             
\begin{proposition}
Let $k = 0$. Then the representations of $\Hh_1$ are
the following:
\begin{itemize}
\item For $\mu,d \in \k, d \neq 0, b = (\mu^p - \mu)^2$,
we have $V_{1,1}^{\mu,d}$, defined as in Proposition \ref{prop3}.
\item For $c, \theta = \pm 1$, we have a $p$-dimensional 
representation $V_{1,6}^{c, \theta}$ with 
basis $\{ v_j, j = 0, 1, \ldots, p-1 \}$, defined by
\begin{eqnarray}
\y v_j &=& j v_j, \quad j = 0, 1, \ldots, p-1; \label{v16first} \\
\s v_0 &=& c v_0; \\
\s v_j &=& - j v_{-j}, \quad j = 1, \ldots, \frac{p-1}{2}; \\
\s v_{-j} &=& - \frac{1}{j} v_j, \quad j = 1, \ldots, \frac{p-1}{2}; \\
\X v_{j} &=& \s v_{-j-1}, \quad j \neq \frac{p-1}{2}; \\
\X v_{\frac{p-1}{2}} &=& \theta \s v_{\frac{p-1}{2}}. \label{v16last}
\end{eqnarray}
\item For $c = \pm 1, a \in \k$, we have a $2p$-dimensional 
representation $V_{1,7}^{c,a}$
with basis $\{ v_j, u_j, j = 0, 1, \ldots, p-1 \}$, defined by
\begin{eqnarray} 
\y v_j &=& j v_j, \quad j = 0, 1, \ldots, p-1; 
\label{v17first} \\
\y u_j &=& j u_j, \quad j = 0, 1, \ldots, p-1; \\
\s v_0 &=& v_0; \\
\s u_0 &=& a v_0 - u_0;
\end{eqnarray}
\begin{eqnarray}
\s v_{-j} &=& - \frac{1}{j} v_j \quad j = 1, 2, \ldots, \frac{p-1}{2}; \\
\s v_j &=& - j v_{-j} \quad j = 1, 2, \ldots, \frac{p-1}{2}; \\
\s u_j &=& \frac{1}{j} u_{-j}, \quad j = 1, 2, \ldots, \frac{p-1}{2}; \\
\s u_{-j} &=& j u_j, \quad j = 1, 2, \ldots, \frac{p-1}{2}; \\
\X v_{j} &=& \s v_{-j-1}, \quad j \neq \frac{p-1}{2}; \\
\X v_{\frac{p-1}{2}} &=& \s u_{\frac{p-1}{2}}; \\
\X u_{j} &=& \s u_{-j-1}, \quad j \neq \frac{p-1}{2}; \\
\X u_{\frac{p-1}{2}} &=& \s v_{\frac{p-1}{2}}. \label{v17last}
\end{eqnarray}
\end{itemize}
$V_{1,1}^{\mu,d}$ and $V_{1,1}^{\mu',d'}$ are isomorphic if
and only if
$$(\mu' - \mu \in \pf \mbox{ {\rm and} } d' = d) \mbox{ {\rm or} }
(\mu' + \mu \in \pf \mbox{ {\rm and} } dd' = \prod_{c \in \pf} 
\left(\frac{k^2}{4} - (\mu+c)^2 \right)).$$
$V_{1,6}^{c,\theta}$ and $V_{1,6}^{c',\theta'}$ are isomorphic if and only if
$\theta = \theta'$.
$V_{1,7}^a$ and $V_{1,7}^{a'}$ are isomorphic if and only if $a = a'$.
Furthermore, representations with different subscripts are never
isomorphic.
\label{prop5}
\end{proposition}

\section{Proof of Propositions \ref{prop1} and \ref{prop2}}

\begin{lemma}[PBW for $\Hh_0$, easy direction]
The elements $$\s^i \X^j \y^l, \quad \quad j,l \in \Z, l \geq 0, i \in 
\{0,1\}$$ span $\Hh_0$ over $\k$.
\label{PBW}
\end{lemma}

\begin{proof}
Given a product of $\X, \y, \s, \X^{-1}$ in any order, one can ensure that 
the
$\y$'s are to the right of all the $\X$'s by using $\y\X = \X\y - k\s\X$ 
repeatedly, and one can also ensure that 
the $\s$'s are to the left of all the $\X$'s and $\y$'s by using $\X\s = 
\s \X^{-1}$ and $\y\s = -k -\s\y$ repeatedly. 
\end{proof}

\begin{lemma}
$\X + \X^{-1}, \y^2$ and $\X\y - \y\X^{-1}$ 
belong to the center 
$\cent{\Hh_0}$ of $\Hh_0$.
\label{lemma1}
\end{lemma}

\begin{proof}
First, let us show that $\y^2 \in \cent{\Hh_0}$. We have
\begin{eqnarray*}
\X \y^2 &=& (\y \X + k \s\X) \y \\
&=& \y\X\y + k\s\X\y \\
&=& \y(\y\X + k\s\X) + k\s(\y\X + k\s\X) \\ 
&=& \y^2 \X + k (\y\s+\s\y) \X + k^2\s^2\X \\
&=& \y^2 \X + -k^2 \X + k^2 \X \\
&=& \y^2 \X;
\end{eqnarray*}
thus $[\X, \y^2] = 0$. We also have
\begin{eqnarray*}
\s \y^2 &=& (-\y\s -k) \y = -\y\s\y - k\y = -\y(-\y\s - k) - k \y
= \y^2 \s;
\end{eqnarray*}
thus $[\s, \y^2] = 0$. It follows that $\y^2 \in \cent{\Hh_0}$.
Next, we show that $\X + \X^{-1} \in \cent{\Hh_0}$.
We have
$$\y (\X + \X^{-1}) = \X\y - k \s\X + \X^{-1} \y + k \X^{-1} \s \\
= (\X + \X^{-1}) \y,$$
and $$\s (\X + \X^{-1}) = \X^{-1} \s + \X \s = (\X + \X^{-1}) \s.$$
Thus $[\y, \X + \X^{-1}] = [\s, \X + \X^{-1}] = 0$, and so $\X + \X^{-1} 
\in \cent{\Hh_0}.$
Finally, we show that $\X\y - \y\X^{-1} \in \cent{\Hh_0}$.
First we note that
$$\y \X - \X^{-1} \y = \y \X + \X \y - (\X + \X^{-1}) \y 
= \y \X + \X \y - \y (\X + \X^{-1}) ) = \X \y - \y \X^{-1},$$
and thus
$$\X (\X \y - \y \X^{-1}) = \X (\y \X - \X^{-1} \y) =
(\X \y - \y \X^{-1}) \X,$$ 
$$\y (\X \y - \y \X^{-1}) = \y (\y \X - \X^{-1} \y) =
\y^2 \X - \y \X^{-1} \y = \X \y^2 - \y \X^{-1} \y
= (\X \y - \y \X^{-1}) \y,$$
and $$\s (\X \y - \y \X^{-1}) = \s (\y \X - \X^{-1} \y) =
- (\y \s + k) \X - \X \s \y = - \y \X^{-1} \s - k \X + \X (\y\s + k)
= (\X \y - \y \X^{-1}) \s.$$
Thus $\X \y - \y \X^{-1} \in \cent{\Hh_0}.$
\end{proof}

\begin{corollary}
$\Hh_0$ is finitely generated as a module over its center.
\label{cor2}
\end{corollary}

\begin{proof}
From Lemmas
\ref{PBW} and \ref{lemma1}, we see that $\Hh_0$ is generated over its 
center by 
$$\s^i \X^j \y^l, \quad \quad i,j,l \in \{0,1\}.$$  
\end{proof}

\begin{corollary}
Every irreducible $\Hh_0$-module is finite-dimensional over $\k$.
\label{cor3}
\end{corollary}

\begin{proof}
Standard.
%Let $V$ be an irreducible $\Hh_0$-module. Then $V$ is cyclic; that is,
%$V = \Hh_0 / J$, where $J$ is a left ideal of $\Hh_0$. By corollary 
%\ref{cor2}, it follows that $V$ is a finitely generated 
%$\cent{\Hh_0}$-module. Hence, there exists a maximal ideal $M \subset 
%\cent{\Hh_0}$ such that $V / MV \neq 0$. But $MV$ is a 
%submodule of the irreducible $\Hh_0$-module $V$, so $MV = 0$. 
%Thus the action of $\Hh_0$ in $V$ factors
%through $\Hh_0/M \Hh_0$, which is a finite-dimensional algebra.
%So $V=(\Hh_0/M \Hh_0)/L$ for some left ideal $L$ of $\Hh_0/M\Hh_0$. 
%Hence $V$ is finite-dimensional.
\end{proof}

Thus Schur's lemma implies that 
central elements of $\Hh_0$ act as scalars in 
any irreducible $\Hh_0$-module. 

From this point, we will use the following notation: the eigenspace of $\y$ with 
eigenvalue $\beta$ will be denoted $V[\beta]$.

\begin{corollary}
Let $V$ be an irreducible $\Hh_0$-module, and let $\beta$ be an 
eigenvalue of $\y$.
Suppose $\beta \neq 0$.
Then, $$V = V[\beta] \oplus V[-\beta],$$ 
and $\dm V[\beta] = \dm V[-\beta] = 1.$
\label{corlast}
\end{corollary}

\begin{proof}
Suppose $V[\beta] \neq 0$, and let $v \in V[\beta]$ be nonzero. From the 
proof of 
corollary \ref{cor2}, we know that $V$ is spanned by
$$\{v, \X v, \s v, \s \X v\}.$$
Now let $w = \s \X v;$ then
$$\y w = \y \s \X v = - \s \y \X v - k \X v 
= - \s \X \y v + k \s^2 \X v - k \X v = - \beta \s \X v = - \beta w.$$ 
Thus, $w \in V[-\beta].$ Clearly, $w \neq 0$, and thus
$V[-\beta] \neq 0$. 
Now let $v' = 2 \beta \X v - k w.$
Then,
$$\y v' = 2 \beta \y \X v + \beta k w = 2 \beta \X \y v - 2 k \beta \s \X 
v + \beta k w = 2 \beta^2 \X v - \beta k w = \beta v'.$$
Hence, $v' \in V[\beta].$ Also, if $w' = k v + 2 \beta \s v,$ then
$$\y w' = k \y v + 2 \beta \y \s v = \beta k v - 2 \beta \s \y v
- 2 \beta k v = - \beta k v - 2 \beta \s \y v = - \beta w',$$
Hence, $w' \in V[-\beta].$
From this it follows that
$$V = V[\beta] \oplus V[-\beta].$$
Now let $\subH_0$ be the subalgebra of $\Hh_0$ generated by $\cent{\Hh_0}$ 
and $2 \beta \X - k \s \X$.
It is clear that $V[\beta] = \subH v$.
Since this is true for all nonzero $v \in V[\beta]$, it follows 
that $V[\beta]$ is an irreducible representation of $\subH_0$. Since 
$\subH_0$ is commutative, we see that $V[\beta]$ is one dimensional. 
The same holds for $V[-\beta]$, and the corollary is proved. 
\end{proof}

\begin{corollary}
Assume $k \neq 0$.
Let $V$ be an irreducible $\Hh_0$-module, and suppose $0$ is an 
eigenvalue of $\y$.
Then, $V = V_{\text{gen}}[0],$ the generalized eigenspace of 
$0$. We also have
$\dm V  = 2$ and $\dm V[0] = 1.$
\label{corlast2}
\end{corollary}

\begin{proof}
Let $v \in V[0]$ be nonzero. From the
proof of corollary \ref{cor2}, we know that $V$ is spanned by
$$\{v, \X v, \s v, \s \X v \}.$$
Let $w = - \s v;$ then,
$$\y w = - \y \s v = \s \y v + k v = k v.$$
Let $v' = \s \X v = \X^{-1} \s v;$ then,
$$\y v' = \y \X^{-1} \s v = \X^{-1} y \s v + k \X^{-1} \s^2 v
= - \X^{-1} \s \y v = 0.$$
Let $w' = - \X v = - \s v';$ then, 
as above, we have $\y w' = v'.$

So we have
$$\y v = 0, \quad \y w = k v, \quad \y v' = 0, \quad \y w' = k v';$$
therefore, $V = V_{\text{gen}}[0]$ and $V[0]$ is spanned by $v$ 
and $v'$.
Now let $\subH_0$ be the subalgebra of $\Hh_0$ generated by $\cent{\Hh_0}$
and $\s \X$.
It is clear that $V[0] = \subH v$.
Since this is true for all nonzero $v \in V[0]$, it follows
that $V[0]$ is an irreducible representation of $\subH_0$. Since
$\subH_0$ is commutative, we see that $V[0]$ is one dimensional.
The corollary follows from this.
\end{proof}

\begin{corollary}
Assume $k = 0$.
Let $V$ be an irreducible $\Hh_0$-module, and suppose $0$ is an
eigenvalue of $\y$.
Then, $$V = V[0],$$ the eigenspace
of $0$. We also have 
$$\dm V = \begin{cases}
1 & \mbox{if } 1 \mbox{ or } -1 \mbox{ is an eigenvalue of } \X \\
2 & \mbox{otherwise}.
\end{cases}$$
\label{corlast3} 
\end{corollary}
                                                                                
\begin{proof}
From the proof of corollary \ref{corlast2}, we see that $\y$ acts on $V$
as the zero operator. Let $\lambda$ be an eigenvalue of $\X$, let 
$\vx[\lambda]$ denote the associated eigenspace and let 
$v \in \vx[\lambda]$ be nonzero. From the
proof of corollary \ref{cor2}, we know that $V$ is spanned by
$\{v, \s v  \}.$ Now 
$$\X \s v = \s \X^{-1} v = \lambda^{-1} \s v,$$
so $\s v \in \vx[\lambda^{-1}].$ Clearly, $\s v \neq 0;$ thus, if
$\lambda \neq \pm 1,$ then $$V = \vx[\lambda] \oplus \vx[\lambda^{-1}] 
\quad \mbox{ and } \dm V = 2.$$
If $\lambda = \pm 1,$ it follows that $\X$ and $\s$ commute as operators 
on $V$; since $V$ is irreducible, this implies that $\dm V = 1.$
\end{proof} 

\begin{proof}[Proof of Proposition \ref{prop1}]
Let $\beta \neq 0$, and let $V$ be a two-dimensional representation of $\Hh_0$ in 
which $V[\beta]$ and $V[-\beta]$ both have dimension 1.
Let $v_0 \in V[\beta]$, $v_1 \in V[-\beta]$ be nonzero.
Let the matrices representing $\s$ and $\X$ with respect
to the basis $\{ v_0, v_1 \}$ be as follows:
$$\s \mapsto
\left( \begin{array}{cc}
\gamma_0 & \delta_0 \\
\gamma_1 & \delta_1
\end{array} \right), \quad
\X \mapsto
\left( \begin{array}{cc}
\theta_0 & \omega_0 \\
\theta_1 & \omega_1
\end{array} \right).
$$
First, we note that $\X$ and $\y$ cannot have a common eigenvector; for 
if $\X w = \gamma w$ and $\y w = \beta' w$, then
$k \s w = \X \y \X^{-1} w - \y w = 0,$ and combining this with $\s^2 = 1$ 
gives $w = 0.$ Hence, by scaling, we can assume that $\omega_0 = 1.$

Now the central element $\X \y - \y \X^{-1}$ acts on $V$ as a scalar.
The matrix representation of $\X \y - \y \X^{-1}$ is
$$\left( \begin{array}{cc}
\frac{\beta}{\det \X} \left(\theta_0 \det \X - \omega_1 \right) 
& -\frac{\beta}{\det \X} (\det \X - 1)  \\
\frac{\beta \theta_1}{\det \X} (\det \X - 1)
& -\frac{\beta}{\det \X} \left(\theta_0 \det \X - \theta_0 \right)
\end{array} \right).$$

Thus, $0 = -\frac{\beta}{\det \X} (\det \X - 1),$ which means that 
$\det \X = 1$.
Hence, 
\begin{eqnarray}
\theta_1 = \theta_0 \omega_1 - 1.
\label{p1e1}
\end{eqnarray}

Using \eqref{rel4}, we see that $\X \y \X^{-1} - \y - k \s = 0.$ 
Using \eqref{p1e1}, 
we see that the matrix 
yrepresentation of $\X \y \X^{-1} - \y - k \s$ is
$$\left( \begin{array}{cc}
2 \beta \theta_0 \omega_1 - 2 \beta - k \gamma_0
& - 2 \beta \theta_0 - k \delta_0 \\
2 \beta \theta_0 \omega_1^2 - 2 \beta \omega_1 - k \gamma_1 
& - 2 \beta \theta_0 \omega_1 + 2 \beta - k \delta_1
\end{array} \right).$$
Hence, 
\begin{eqnarray}
\gamma_0 &=& \frac{2 \beta}{k} (\theta_0 \omega_1 - 1), \label{p1e2} \\
\gamma_1 &=& \frac{2 \beta}{k} \left(\theta_0 \omega_1^2 - \omega_1 
\right), 
\label{p1e3} \\
\delta_0 &=& - \frac{2 \beta}{k} \theta_0 \label{p1e4} \\
\delta_1 &=& \frac{2 \beta}{k} (1 - \theta_0 \omega_1). \label{p1e5} 
\end{eqnarray}

Using \eqref{rel2}, we see that $\s^2 = 1$. 
Using \eqref{p1e2}--\eqref{p1e5}, we see 
that the matrix representation of $\s^2$ is
$$\left( \begin{array}{cc}
\frac{4 \beta^2}{k^2} (1 - \theta_0 \omega_1) & 0 \\
0 & \frac{4 \beta^2}{k^2} (1 - \theta_0 \omega_1)
\end{array} \right).$$
Thus,
\begin{equation} 
\omega_1 = \frac{1}{\theta_0} \left(1 - \frac{k^2}{4 \beta^2} \right). 
\label{p1e6}
\end{equation}

Using \eqref{p1e1}--\eqref{p1e6}, we see that $V$ is isomorphic to
$V_{0,1}^{\beta,\theta_0}.$
Furthermore, it is easy to see that
for all $a, \beta \in \k \nozero,$ $V_{0,1}^{\beta, a}$ is a 
representation of $\Hh_0$; furthermore, each eigenvector of $\y$ clearly 
generates $\Hh_0$, and thus $V_{0,1}^{\beta,\theta_0}$ is irreducible.

Now the eigenvalues of $\y$ in $V_{0,1}^{\beta, a}$ 
are $\beta$ and $-\beta$, and $2 \beta(\X \y - \y \X^{-1}) - 2 \beta^2(\X 
+ \X^{-1})$ acts on $V_{0,1}^{\beta, a}$ as $\frac{k^2 - 4 \beta^2}{a} 
\Id,$ while 
for $\beta = \frac{k}{2},$ $\X + \X^{-1}$ acts on $V_{0,1}^{\beta, a}$ as 
$a \Id$. From this it follows that $V_{0,1}^{\beta, a}$ and 
$V_{0,1}^{\beta', a'}$ are 
isomorphic if and only if $\beta' = \beta, a' = a$ or $\beta' = -\beta, a' 
= \frac{4 \beta^2 - k^2}{4 a \beta^2}$.

Now let $V$ be a two-dimensional representation of $\Hh_0$ in which $V[0]$ 
has dimension 1 and $\vgen[0]$ has dimension 2. Let $v_0, v_1 \in V$ be 
nonzero elements such 
that $\y v_0 = 0, \y v_1 = v_0$. 
Let the matrices representing $\s$ and $\X$ with respect
to the basis $\{ v_0, v_1 \}$ be as follows:
$$\s \mapsto
\left( \begin{array}{cc}
\gamma_0 & \delta_0 \\
\gamma_1 & \delta_1
\end{array} \right), \quad
\X \mapsto
\left( \begin{array}{cc}
\theta_0 & \omega_0 \\
\theta_1 & \omega_1
\end{array} \right).
$$
Now the matrix representation of $\s \y + \y \s + k$ is
$$
\left( \begin{array}{cc}
\gamma_1 + k & \gamma_0 + \delta_1 \\
0 & \gamma_1 + k
\end{array} \right).$$
\eqref{rel3} thus implies that $\gamma_1 = -k$ and $\gamma_0 = -\delta_1.$
Scaling, we may assume that $\gamma_0 = 1$.
Next, we note that $\s^2 - 1$ acts on $V$ as $-\delta_0 k \Id;$ 
\eqref{rel2} thus implies that $\delta_0 = 0.$
We then see that the matrix representation of $\X\y - \y\X - k\s\X$ is
$$
\left( \begin{array}{cc}
-\theta_1 - k \theta_0 & \theta_0 - \omega_1 - k \omega_0 \\
k ( k \theta_0 + \theta_1 ) & \theta_1 + k^2 \omega_0 + k \omega_1
\end{array} \right).$$
\eqref{rel4} thus implies that $\theta_1 = -k \theta_0, \theta_0 = 
\omega_1 + k \omega_0.$
Finally, the matrix representation of $\X \s \X - \s$ is
$$
\left( \begin{array}{cc}
\theta_0^2 - 1 & 0 \\
- k ( \theta_0^2 - 1 ) & - \theta_0^2 + 1 
\end{array} \right).$$
\eqref{rel4} thus implies that $\theta_0 = \pm 1$.

Thus $V$ is isomorphic to $V_{0,2}^{\theta_0, \omega_0}.$
It is easy to see that $V_{0,2}^{a,b}$ is indeed a representation 
of $\Hh_0$; furthermore, each eigenvector of $\y$ clearly generates $\Hh_0$, and
thus $V_{0,2}^{a,b}$ is irreducible.
Now $\X + \X^{-1}$ acts on $V_{0,2}^{a,b}$ as $(2 a
- kb) \Id,$ while $\X\y - \y\X^{-1}$ acts as $- ak \Id.$ 
Therefore, $V_{0,2}^{a,b}$ and $V_{0,2}^{a', b'}$ are 
isomorphic if and only if $a' = a, b' = b.$ 
\end{proof}

\begin{proof}[Proof of Proposition \ref{prop2}]
Let $\beta \neq 0$, and let $V$ be a two-dimensional representation of 
$\Hh_0$ in which $V[\beta]$ and $V[-\beta]$ both have dimension 1.
Let $v_0 \in V[\beta]$, $v_1 \in V[-\beta]$ be nonzero.
Let the matrices representing $\s$ and $\X$ with respect
to the basis $\{ v_0, v_1 \}$ be as follows:
$$\s \mapsto
\left( \begin{array}{cc}
\gamma_0 & \delta_0 \\
\gamma_1 & \delta_1
\end{array} \right), \quad
\X \mapsto
\left( \begin{array}{cc}
\theta_0 & \omega_0 \\
\theta_1 & \omega_1
\end{array} \right).
$$
First, we note thet $\X$ and $\y$ commute, so they must have a 
common eigenvector; for the moment, let us assume that $\omega_0 = 0$. 
                                                                                
Now the central element $\X \y - \y \X^{-1}$ acts on $V$ as a scalar.
The matrix representation of $\X \y - \y \X^{-1}$ is
$$\left( \begin{array}{cc}
\frac{\beta}{\det \X} (\theta_0 \det \X - \omega_1) & 0  \\
\frac{\beta \theta_1}{\det \X} (\det \X - 1)
& -\frac{\beta}{\det \X} (\omega_1 \det \X - \theta_0)
\end{array} \right).$$
                                                                                
Thus, $0 = \frac{\beta \theta_1}{\det \X} (\det \X - 1),$ which means that
$\det \X = 1$.
Hence,
\begin{eqnarray}
\omega_1 = \frac{1}{\theta_0}.
\label{p2e1}
\end{eqnarray}
                                                                                
Using \eqref{rel4}, we see that $\X \y \X^{-1} - \y = 0.$
Using \eqref{p2e1},
we see that the matrix
representation of $\X \y \X^{-1} - \y$ is
$$\left( \begin{array}{cc}
0 & 0 \\
\frac{2 \beta}{\theta_0} \theta_1
& 0
\end{array} \right).$$
Hence, $\theta_1 = 0.$ (If we had assumed earlier that $\theta_1=0$, here 
we would get $\omega_0=0.$)

Using \eqref{rel1}, we see that $\y \s + \s \y = 0.$
Using \eqref{p2e1},
we see that the matrix
representation of $\y \s + \s \y$ is
$$\left( \begin{array}{cc}
2 \beta \gamma_0 & 0 \\
0 & - 2 \beta \delta_1
\end{array} \right).$$

Then we must have $\gamma_0 = 
\delta_1 = 0,$ and thus $\s^2$ acts on $V$ as $\delta_0 \gamma_1 
\Id.$ Using \eqref{rel2}, we see that $\s^2 = 1;$ thus $\delta_0 = 
\frac{1}{\gamma_1}.$ By scaling, we may assume that $\gamma_1 = 1,$ and we
see that
$V = V_{0,3}^{\beta,\theta_0}.$ It is clear that $V_{0,3}^{\beta, a}$ 
is an irreducible representation of 
$\Hh_0$ and that $V_{0,3}^{\beta, a}$ and 
$V_{0,3}^{\beta', a'}$ are 
isomorphic if and 
only if $\beta' = \beta, a' = a$ or
$\beta' = -\beta, a' = \frac{1}{a}.$

Now let $V$ be a two-dimensional 
representation of $\Hh_0$ in which $\y$ acts as zero and $\X$ has 
eigenvalues $\lambda$ and $\lambda^{-1}$, with $\lambda \neq \pm 1$.
Let $v_0 \in \vx[\lambda]$, $v_1 \in \vx[\lambda^{-1}]$ be nonzero.
Let the matrix representing $\s$ with respect
to the basis $\{ v_0, v_1 \}$ be as follows:
$$\s \mapsto
\left( \begin{array}{cc}
\gamma_0 & \delta_0 \\
\gamma_1 & \delta_1
\end{array} \right). 
$$

Using \eqref{rel1}, we see that $\X \s - \s \X^{-1} = 0.$ But the matrix 
representation of $\X \s - \s \X^{-1}$ is 
$$
\left( \begin{array}{cc}
\gamma_0 \frac{\lambda^2 - 1}{\lambda} & 0 \\
0 & - \delta_1 \frac{\lambda^2 - 1}{\lambda}
\end{array} \right).
$$
Since $\lambda \neq \pm 1$, we see that $\gamma_0 = \delta_1 = 0.$ Thus 
$\s^2$ acts on $V$ as $\gamma_1 \delta_0 \Id$. Using \eqref{rel2}, we see 
that $\s^2 = 1.$ Hence, $\delta_0 = \frac{1}{\gamma_1}.$ By scaling, 
we may assume that $\gamma_1 = \delta_0 = 1.$ Thus $V$ must
be isomorphic to $V_{0,4}^{\lambda}.$ It is clear that $V_{0,4}^a$ 
is an irreducible representation of $\Hh_0$ and that $V_{0,4}^a$ and 
$V_{0,4}^{a'}$ 
are isomorphic if and only if $a' = a$ or
$a' = \frac{1}{a}.$

Finally, the classification of one-dimensional representations of $\Hh_0$ 
is trivial.
\end{proof}

\section{Proof of Propositions \ref{prop3}, \ref{prop4} and \ref{prop5}}

\begin{lemma}[PBW for $\Hh_1$, easy direction]
The elements $$\s^i \X^j \y^l, \quad \quad j,l \in \Z, l \geq 
0, i \in
\{0,1\}$$ span $\Hh_1$ over $\k$.
\label{1PBW}
\end{lemma}

\begin{proof}
Similar to the proof of lemma \ref{PBW}.
\end{proof}

\begin{lemma}
$\X^{p} + \X^{-p}$ and $(\y^p - \y)^2$
belong to the center
$\cent{\Hh_1}$ of $\Hh_1$.
\label{1lemma1}
\end{lemma}

\begin{proof}
First, let us show that $\X^{p} + \X^{-p} \in \cent{\Hh_1}$. We have
\begin{eqnarray*}
\y (\X^p + \X^{-p})  
&=& \y (\X + \X^{-1})^p \\
&=& (\X \y + \X - k \s \X + \X^{-1} \y - \X^{-1} + k \X^{-1} \s) (\X 
+ \X^{-1})^{p-1} \\
&=& (\X + \X^{-1}) \y (\X^{p-1} + \X^{-p+1}) + (\X - \X^{-1}) (\X
+ \X^{-1})^{p-1} \\
&=& \cdots \\
&=& (\X + \X^{-1})^p \y + p (\X - \X^{-1}) (\X + \X^{-1})^{p-1} \\
&=& (\X^p + \X^{-p}) \y;
\end{eqnarray*} 
thus $[\y, \X^p + \X^{-p}] = 0$. We also have
$$\s \X^p = \X^{-p} \s, \quad \s \X^{-p} = \X^p \s;$$ 
thus $[\s, \X^p + \X^{-p}] = 0$. It follows that $\X^p + \X^{-p} \in 
\cent{\Hh_1}$.
Next, we show that $(\y^p - \y)^2 \in \cent{\Hh_1}$.
We have
\begin{eqnarray*}
\X (\y+1)^2
&=& \X \y^2 + 2 \X \y + \X \\
&=& (\y \X - \X + k \s \X) \y + 2 \X \y + \X \\
&=& \y (\y \X - \X + k \s \X) - \X \y + k \s (\y \X - \X + k \s \X)
+ 2 \X \y + \X \\
&=& \y^2 \X + k (\y \s + \s \y + k) \X + (\X \y - \y \X + \X - k \s \X) \\
&=& \y^2 \X. 
\end{eqnarray*}
So $\X (\y+1)^2 = \y^2 \X,$ and thus $\X g(\y+1) = g(\y) \X$ for all even 
polynomials $g.$ In particular, $[\X, (\y^p - \y)^2] = 0.$
Furthermore, we have $[\s, \y^2] = 0$ (for the same reason as in the case 
$t=0).$ It follows that $(\y^p - \y)^2 \in \cent{\Hh_1}.$
\end{proof}

\begin{corollary}
$\Hh_1$ is finitely generated as a module over its center.
\label{1cor2}
\end{corollary}
                                                                                
\begin{proof}
From Lemmas
\ref{1PBW} and \ref{1lemma1}, we see that $\Hh_1$ is generated over 
its center by
$$\s^i \X^j \y^l, \quad \quad i \in \{0,1\}, j \in \{-p+1, -p+2, 
\ldots, p-1, p \}, l \in
\{0, 1, \ldots, 2p-1 \},$$
\end{proof}
                                                                   
\begin{corollary}
Every irreducible $\Hh_1$-module is finite-dimensional over $\k$.
\label{1cor3}
\end{corollary}

\begin{proof}
Standard.
%Similar to the proof of corollary \ref{cor3}.
\end{proof}

Consider the 
following elements of $\Hh_1:$
$$\iA \defeq \s\X, \quad \iB = \s\y + \frac{k}{2}.$$
These elements were introduced by Cherednik in \cite{Ch4} and are called 
{\em intertwiners}.
We note that $\iB$ is also equal to $-\y\s - \frac{k}{2}.$ 

\begin{lemma}
$$\iA^2 = 1, \quad \iB^2 = - \y^2 + \frac{k^2}{4}.$$
\end{lemma}

\begin{proof}
We have
$$\iA^2 = \s\X\s\X = \s \s \X^{-1} \X = 1$$
and
$$\iB^2 = (\s\y + k) \s\y + \frac{k^2}{4}
= -\y\s \s\y + \frac{k^2}{4} = -\y^2 + \frac{k^2}{4}.$$
\end{proof}

\begin{lemma}
$$\iA \y = (-\y -1) \iA, \quad \iB \y = - \y \iB.$$
\end{lemma}

\begin{proof}
We have
$$\iA \y = \s \X \y = \s (\y \X - \X + k \s \X)
= - \y \s \X - k \X - \s \X + k \X = (-\y - 1) \iA,$$
and
$$\iB \y = - \y \s \y - \frac{k}{2} \y =   
\y^2 \s + \frac{k}{2} \y = - \y \iB.$$
\end{proof}

\begin{corollary}
Let $V$ be a representation of $\Hh_1$.
% and let $\beta$ be the  
%eigenspace of $\y$ with eigenvalue $\beta$.
Then $$\iA : V[\beta] \to V[\-\beta-1]$$ is an isomorphism and $$\iB : 
V[\beta] \to V[\-\beta]$$ is a homomorphism. $\iB$ is an isomorphism if 
and only if $\beta \neq \pm \frac{k}{2}.$ The same result holds for 
generalized eigenspaces.
\label{1cor4}
\end{corollary}

\begin{lemma}
Let $V$ be an irreducible representation of $\Hh_1$ on which the central 
element $(\y^p - \y)^2$ acts as $b \neq 0$. Then
$$V \supset \oplus_{c \in \pf} \left( V[\mu + c] \oplus V[-\mu+c] \right),$$ 
where $\mu$ is a root of the equation $(\mu^p - \mu)^2 = b.$
Each eigenspace has dimension $1$.
%, and thus $\dm V \geq 2p.$
\label{1lem6}
\end{lemma}

\begin{proof}
Let $\mu$ be an eigenvalue of $\y$, and let $v \in V[\mu].$ Note that $b v 
= (\y^p - \y)^2 v = (\mu^p - \mu)^2 v,$ so we have $(\mu^p - \mu)^2 = b.$
By corollary \ref{1cor4}, we have the following homomorphisms:
\begin{equation}
V[\mu] \stackrel{\iA}{\to} V[-\mu-1] \stackrel{\iB}{\to} 
V[\mu+1] \stackrel{\iA}{\to} V[-\mu-2] \stackrel{\iB}{\to} V[\mu+2]
\stackrel{\iA}{\to} \cdots \stackrel{\iA}{\to} V[-\mu-p+1] 
\stackrel{\iB}{\to} V[\mu+p-1] \stackrel{\iA}{\to} V[-\mu]. 
\label{eqn91}
\end{equation}
%Note that if $\mu'$ is one of the eigenvalues in \eqref{eqn91}, 
%then $(\mu'^p - \mu')^2 = b$; hence $\mu' \notin \pf.$ This means that
%the eigenvalues in \eqref{eqn91} are all distinct and that none of them is 
%zero; thus lemma \ref{1lem5} applies, and
%$$V = \oplus_{c \in \pf} \left( V[\mu + c] \oplus V[-\mu+c] \right).$$
If none of the eigenvalues in \eqref{eqn91} is equal to $\frac{k}{2},$ 
then all of the homomorphisms in \eqref{eqn91} are 
isomorphisms (by corollary \ref{1cor4}). 
Otherwise, we may assume without loss of generality that $\mu = 
\frac{k}{2},$ and once again all of the homomorphisms in \eqref{eqn91} 
are isomorphisms.
Thus, $\dm V \geq 2p \dm V[\mu]$.
Now the dimension of the algebra $\Hh_1 / (\X^p + \X^{-p} = a, 
(\y^p - \y)^2 = b)$ acting irreducibly on $V$ is at most $8 p^2$ (see the 
proof of corollary \ref{1cor2}). Hence, $4 p^2 (\dm V[\mu])^2 = (\dm V)^2 
\leq 8 p^2,$ which implies that $\dm V[\mu] = 1.$ The result follows.
\end{proof}

\begin{lemma}
Suppose $k \notin \pf.$
Let $V$ be an irreducible representation of $\Hh_1$ on which the central
element $(\y^p - \y)^2$ acts as $0$. Then
each generalized eigenspace $\vgen[c], c \in \pf$
has dimension $2$.
%, and thus $\dm V \geq 2p.$
\label{1lem8}
\end{lemma}

\begin{proof}
Let $\mu$ be an eigenvalue of $\y$, and let $v \in V[\mu].$ Note that $0
= (\y^p - \y)^2 v = (\mu^p - \mu)^2 v,$ so we have $\mu^p - \mu = 0.$
Hence, $\mu \in \pf.$ 
By corollary \ref{1cor4}, we have the following homomorphisms:
\begin{equation} 
\vgen[\mu] \stackrel{\iB \iA}{\to} \vgen[\mu+1] 
\stackrel{\iB \iA}{\to} \vgen[\mu+2]
\stackrel{\iB \iA}{\to} \cdots 
\stackrel{\iB \iA}{\to} \vgen[\mu-1]
\stackrel{\iB \iA}{\to} \vgen[\mu].
\label{eqn92}
\end{equation}
Since $k \notin \pf,$ none of the eigenvalues in \eqref{eqn92} is equal to
$\pm \frac{k}{2}.$ By corollary \ref{1cor4}, all of the homomorphisms in 
\eqref{eqn92} are isomorphisms, and the eigenvalues of $\y$ are precisely 
the elements of $\pf.$
Thus, $\dm V = p \dm \vgen[0]$.
Now the dimension of the algebra $\Hh_1 / (\X^p + \X^{-p} = a,  
(\y^p - \y)^2 = 0)$ acting irreducibly on $V$ is at most $8 p^2$ (see the
proof of corollary \ref{1cor2}). Hence, $p^2 (\dm \vgen[0])^2 = (\dm V)^2
\leq 8 p^2,$ which implies that $\dm \vgen[0] \leq 2.$ 
Now let $v \in V[0];$ then $$\y \s v = - \s \y v - k v = - k v.$$
Since $k \neq 0,$ we conclude that $\s v \in \vgen[0] \setminus V[0].$
Therefore, $\dm \vgen[0] = 2,$ and the result follows.
\end{proof}

\begin{proof}[Proof of Proposition \ref{prop3}]
It is easy to show that if $\mu$ and $d$ satisfy the conditions in the statement of 
Proposition \ref{prop3}, then $V_{1,1}^{\mu,d}$ is a representation of $\Hh_1$.
Furthermore, if $v \in V_{1,1}^{\mu,d}$ is an eigenvector of $\y$, we 
see that we can generate all of $V_{1,1}^{\mu,d}$ by applying $\iA$ and $\iB$. 
This implies that $V_{1,1}^{\mu,d}$ is actually an irreducible representation of 
$\Hh_1$. The same can be said of $V_{1,2}^{\theta}.$

Let $V$ be an irreducible representation of $\Hh_1$, and suppose that 
$(\y^p - \y)^2$ acts on $V$ as $b \neq 0$.
For the moment, let us assume that $\pm \frac{k}{2}$ are not eigenvalues 
of $\y$. Let $v_{\mu}$ be an eigenvector of $\y$ with eigenvalue $\mu$, 
and let
\begin{eqnarray*}
v_{\mu+j} &=& (\iB \iA)^j v_{\mu}, \quad j = 1, 2, \ldots, p-1 \\
v_{-\mu+j} &=& \iA (\iB \iA)^{j-1} v_{\mu}, \quad j = 1, 2, \ldots, p.
\end{eqnarray*}
Note that $\iB v_{-\mu} \in V[\mu]$; using lemma \ref{1lem6}, we see that
$\iB v_{-\mu} = d v_{\mu},$ where $d \in \k$ is nonzero.
From this we can deduce that \eqref{v11first}--\eqref{v11last} are
satisfied. Thus $V_{1,1}^{\mu,d} \subset V.$
By irreducibility of $V$, it follows that $V = V_{1,1}^{\mu,d}.$
Now we note that $(\iB \iA)^p$ acts as $d \Id$ on $V[\mu], V[\mu+1],
\ldots, V[\mu+p-1]$ and as $\frac{1}{d} \prod_{c \in \pf} \left( 
\frac{k^2}{4} - (\mu+c)^2 \right)$ on
$V[-\mu], V[-\mu+1], \ldots, V[-\mu+p-1].$ From this we can deduce that
$V_{1,1}^{\mu,d}$ and $V_{1,1}^{\mu',d'}$ are isomorphic if and only if
($\mu' - \mu \in \pf$ and $d' = d$) or
($\mu' + \mu \in \pf$ and $dd' = \prod_{c \in \pf} \left( 
\frac{k^2}{4} - (\mu+c)^2 \right)$).
Now, if $\pm \frac{k}{2}$ are eigenvalues of $\y$, then $\iB^2$ acts 
as zero on $V[\pm \frac{k}{2}]$, so either $\iB$ acts as zero on 
$V[\frac{k}{2}]$, in
which case we can use the above argument with $\mu = \frac{k}{2}$, or
$\iB$ acts as zero on $V[\frac{k}{2}]$, in
which case we can use the above argument with $\mu = -\frac{k}{2}$. 
%Finally, if $v \in V_{1,1}^{\mu,d}$ is an eigenvector of $\y$, we see that we 
%can generate all of $V_{1,1}^{\mu,d}$ by applying $\iA$ and $\iB$. 
%This implies that $V_{1,1}^{\mu,d}$

Second, suppose that $(\y^p - \y)^2$ acts on $V$ as $0$.
Let $v_0$ be an eigenvector of $\y$ with eigenvalue $0$,
and let
\begin{eqnarray*}
v_j &=& (\iB \iA)^j v_0, \quad j = 1,2, \ldots, \frac{p-1}{2}; \\
v_{-j} &=& - \iA (\iB \iA)^{j-1} v_0, 
\quad j = 1, 2, \ldots, \frac{p-1}{2}; \\
w_0 &=& - \frac{1}{k} \s v_0; \\
w_j &=& (\iB \iA)^j w_0, \quad j = 1, 2, \ldots, \frac{p-1}{2}; \\
w_{-j} &=& \iA (\iB \iA)^{j-1} w_0, 
\quad j = 1, 2, \ldots, \frac{p-1}{2}.
\end{eqnarray*}
Since $\iA$ maps eigenspaces to eigenspaces and generalized eigenspaces to 
generalized eigenspaces, and since $\iA^2=1$, we can use lemma 
\ref{1lem8} to conclude that
$$\iA v_{\frac{p-1}{2}} = \theta v_{\frac{p-1}{2}}, \quad
\iA w_{\frac{p-1}{2}} = \omega w_{\frac{p-1}{2}}, \quad
\theta, \omega = \pm 1.$$
From this we can deduce that \eqref{v12first}--\eqref{v12-2ndlast} are
satisfied, as well as
\begin{equation}
\X w_{\frac{p-1}{2}} = \omega \s w_{\frac{p-1}{2}} \label{expr88}
\end{equation}
Now we know from \eqref{rel4} that 
\begin{equation}
(\X \y - \y \X + \X - k \s \X) w_{\frac{p-1}{2}} \label{expr99}
\end{equation}
must be zero.
Using \eqref{v12first}--\eqref{v12-2ndlast} and \eqref{expr88}, we see 
that the coefficient of $v_{\frac{p-1}{2}}$ in $\X \y w_{\frac{p-1}{2}}, - 
\y \X w_{\frac{p-1}{2}}, \X w_{\frac{p-1}{2}}$ and $- k \s \X 
w_{\frac{p-1}{2}}$ are respectively $k (\theta - \omega), 0, 2 k \omega$ 
and $0$. From \eqref{expr99}, we get $k (\theta + \omega) = 0,$ which 
implies $\omega = -\theta$. Hence \eqref{v12last} is also satisfied. 
Thus $V^{1,2}_{\theta} \subset V.$
By irreducibility of $V$, it follows that $V = V_{1,2}^{\theta}.$ Since
$\iA$ acts on $V_{1,2}^{\theta}[\frac{p-1}{2}]$ through multiplication by 
$\theta$, it is clear that $V_{1,2}^{\theta}$ and $V_{1,2}^{\theta'}$ are 
isomorphic if and only if $\theta=\theta'$.
\end{proof}
                                                                               
\begin{lemma}
Let $k$ be an even integer with $2 \leq k \leq p-1$. Let $V \neq 0$ be an 
irreducible representation of $\Hh_1$ on which $(\y^p - \y)^2$ acts as zero.
Then $V \left[\frac{k}{2} \right] \neq 0$.
\end{lemma}

\begin{proof}
From corollary \ref{1cor4}, we have isomorphisms
\begin{equation}
V[0] \stackrel{\iA}{\to} V[-1] \stackrel{\iB}{\to}
V[1] \stackrel{\iA}{\to} V[-2] \stackrel{\iB}{\to} V[2]
\stackrel{\iA}{\to} \cdots \stackrel{\iA}{\to} V \left[-\frac{k}{2} \right]
\label{isom66}
\end{equation}
and
$$V \left[\frac{k}{2} \right] \stackrel{\iA}{\to} V \left[-\frac{k}{2}-1 \right] 
\stackrel{\iB}{\to} V \left[ \frac{k}{2} +1 \right] \stackrel{\iA}{\to} \cdots 
\stackrel{\iA}{\to} V \left[ \frac{p-1}{2} \right].$$
Let us assume that $V \left[\frac{k}{2} \right] = 0$. Since $V \neq 0$, we must 
have $V[0] \neq 0$.
Let $u \in V[0]$ be nonzero; then $\y \s u = - \s \y u - k u = - k u.$ This means 
$\vgen[0] \setminus V[0]$ is nonempty. From the isomorphisms \eqref{isom66},
we know that there exist $v, w \in V$ such that $\y v = - \frac{k}{2} v, \y w = - 
\frac{k}{2} w + v$. Since $V \left[\frac{k}{2} \right] = 0,$ we see that $\iB v = 
\iB w = 0$. Thus, 
$$- \frac{k}{2} v = \s \y v = - \frac{k}{2} \s v \impl \s v = v;$$
$$- \frac{k}{2} w = \s \y w = - \frac{k}{2} \s w + \s v 
= - \frac{k}{2} \s w + v \impl \s w = \frac{2}{k} v + w.$$
Hence $$w = \s^2 w = \frac{2}{k} \s v + \s w = \frac{4}{k} v + w,$$ which is 
impossible. Therefore, $V \left[\frac{k}{2} \right] \neq 0$.
\end{proof}

\begin{lemma}
Let $V$ be an irreducible representation of $\Hh_1$, 
and suppose that $\iB : V[\frac{k}{2}] 
\to V[-\frac{k}{2}]$ is zero but $V[\frac{k}{2}] \neq 0$. Then
$$V = V[\frac{k}{2}] \oplus V[\frac{k}{2}+1] \oplus \cdots \oplus 
V[-\frac{k}{2}-1].$$
\label{lemma50}
\end{lemma}

\begin{proof}
Let $$W = V[\frac{k}{2}] \oplus V[\frac{k}{2}+1] \oplus \cdots \oplus
V[-\frac{k}{2}-1].$$ Since $V$ is irreducible, it is enough to show that $W$ 
is a subrepresentation. Since $\iB$ acts as zero on $V[\frac{k}{2}],$ we see 
that $W$ is closed under the action of $\iA$ and $\iB$. Clearly, $\y W 
\subset W$. For $v \in V[\beta], \beta \neq 0$, we have
$$\iB v = \left( \s \y + \frac{k}{2} \right) v
= \beta \s v + \frac{k}{2} v; \mbox{ so } \s v = \beta^{-1} \left( \iB - 
\frac{k}{2} \right) v.$$
So $\s W \subset W$. Finally, $\X = \s \iA,$ so $\X W \subset W$, and the 
proof is complete.
\end{proof}

\begin{proof}[Proof of Proposition \ref{prop4}]
As in the proof of Proposition \ref{prop3}, we see that $V_{1,1}^{\mu,d}$ is an 
irreducible representation of $\Hh_1$ whenever $\mu$ and $d$ satisfy the conditions 
in the statement of Proposition \ref{prop4}. Similarly, for 
$\theta = \pm 1,$ $V_{1,3}^{\theta}$ is an irreducible representation of $\Hh_1$
and for $c \in \k$, $V_{1,4}^c$ and $V_{1,5}^c$ are irreducible
representations of $\Hh_1$.

Now, let $V$ be an irreducible representation of $\Hh_1$, 
and suppose that that $(\y^p - \y)^2$ acts on $V$ as $b$.
Here $k \in \pf$, so if $b \neq 0$, then $\frac{k}{2}$ is not a root of 
$f(y) = (y^p - y)^2 - b$, so the argument in the proof of Proposition 
\ref{prop3} applies for $V_{1,1}^{\mu, d}$. We will now assume that $b = 0$.

Now suppose $\iB$ acts as zero on $V[\frac{k}{2}].$
Let $\subH_1$ be the subalgebra of $\Hh_1$ generated by $\iA$ and $\iB$.
From the proof of Lemma \ref{lemma50}, we see that $V = \subH_1 v$ for any 
eigenvector of $v$ of $\y$. Since $\iA^2 = 1$ and $\iB^2$ acts as a scalar 
on each eigenspace, it follows that $V$ is spanned by
$$v, \iA v, \iB \iA v, \iA \iB \iA v, \ldots, \iB v, \iA \iB v, \iB \iA 
\iB v, \ldots$$ for each eigenvector $v$ of $\y$. If the eigenvalue of $v$ 
is $\frac{k}{2},$ we conclude, from corollary \ref{1cor4} and from the fact 
that $\iB$ acts as zero on $V[\frac{k}{2}]$, that 
\begin{equation}
V[\frac{k}{2}] \mbox{ is spanned by } v, \iA (\iB \iA)^{p-k-1} v.
\label{eqn100}
\end{equation}
Let $\subsubH_1$ be the subalgebra of $\subH_1$ generated by $\iA 
(\iB \iA)^{p-k-1}$. Clearly, $\subsubH_1$ is commutative, and since 
\eqref{eqn100} 
holds for all nonzero $v \in V[\frac{k}{2}],$ we conclude that 
$V[\frac{k}{2}]$ is an irreducible representation of $\subsubH$. By 
Schur's Lemma, $\iA (\iB \iA)^{p-k-1}$ acts on $V[\frac{k}{2}]$ as a 
scalar, and 
thus $\dm V[\frac{k}{2}] = 1.$ We then let $v_{\frac{k}{2}} \in 
V[\frac{k}{2}]$ be nonzero, and let
\begin{eqnarray*}
v_{-\frac{k}{2}-j} = \iA (\iB \iA)^{j-1} v, \quad j = 1, 2, \ldots, 
\frac{p-1}{2} - \frac{k}{2}; \\
v_{\frac{k}{2}+j} = (\iB \iA)^j v, \quad j = 1, 2, \ldots,
\frac{p-1}{2} - \frac{k}{2}.
\end{eqnarray*}
Now $\iA v_{\frac{p-1}{2}} = \theta v_{\frac{p-1}{2}}$ for some $\theta,$ 
and $\iA^2 = 1$ implies $\theta = \pm 1.$
From the above information, we can deduce that 
\eqref{v13first}--\eqref{v13last} are satisfied, and thus 
$V_{1,3}^{\theta} \subset V$. 
By irreducibility of $V$, we conclude that $V = V_{1,3}^{\theta}$.
Since $\iA$ acts on $V[\frac{p-1}{2}]$ through multiplication by $\theta$, 
it follows that $V_{1,3}^{\theta}$ and $V_{1,3}^{\theta'}$ are isomorphic 
if and only if $\theta = \theta'$.

Now suppose $\iB$ does not act as zero on $V[\frac{k}{2}].$
Let $v_0$ be an eigenvector of $\y$ with eigenvalue $0$, and let $w_0 = 
-\frac{1}{k} \s v_0.$ Then $\y w_0 = v_0.$
For $j = 0, \ldots, \frac{k}{2} - 1,$ let
\begin{eqnarray*}
v_j &=&  (\iB \iA)^j v_0 \\
w_j &=&  (\iB \iA)^j w_0 \\
v_{-j-1} &=&  - \iA (\iB \iA)^{j} v_0 \\
w_{-j-1} &=&  \iA (\iB \iA)^{j} w_0.
\end{eqnarray*}
It is easy to check that for each $j$, $v_j$ and $v_{-j-1}$ are 
eigenvectors of 
$\y$ with eigenvalue $j, -j-1$ respectively, and $\y w_j = j w_j + v_j, \y 
w_{-j-1} = (-j-1) w_{-j-1} + v_{-j-1}.$ 
Now $\iB v_{-\frac{k}{2}} = 0,$ and 
$$\y \iB w_{-\frac{k}{2}} = - \iB \y w_{-\frac{k}{2}} =
\frac{k}{2} \iB w_{-\frac{k}{2}} - \iB v_{-\frac{k}{2}}
= \frac{k}{2} \iB w_{-\frac{k}{2}}.$$
Thus, $w_{-\frac{k}{2}} \in V[\frac{k}{2}].$ Let us write $v_{\frac{k}{2}}
= \iB w_{-\frac{k}{2}};$
then,
\begin{eqnarray*}
\s w_{-\frac{k}{2}} &=& - \frac{2}{k} v_{\frac{k}{2}} + \frac{2}{k} 
v_{-\frac{k}{2}} + w_{-\frac{k}{2}} \\
\mbox{and } w_{-\frac{k}{2}} = \s^2 w_{-\frac{k}{2}} &=&
- \frac{2}{k} \s v_{\frac{k}{2}} + \frac{2}{k} v_{-\frac{k}{2}} + \s 
w_{-\frac{k}{2}} \\
&=& -\frac{2}{k} \s v_{\frac{k}{2}} + \frac{4}{k} v_{-\frac{k}{2}} - 
\frac{2}{k} z + w_{-\frac{k}{2}} \\
\impl 0 &=& - \s v_{\frac{k}{2}} + 2 v_{-\frac{k}{2}} - v_{\frac{k}{2}}.
\end{eqnarray*}
In particular, this implies that $v_{\frac{k}{2}} \neq 0.$
For $j = \frac{k}{2} + 1 , \ldots, \frac{p-1}{2},$ let
\begin{eqnarray*}
v_{-j} &=&  - \iA (\iB \iA)^{j - \frac{k}{2} - 1} v_{\frac{k}{2}} \\
v_{j} &=& (\iB \iA)^{j - \frac{k}{2}} v_{\frac{k}{2}}.
\end{eqnarray*}
Now there are two cases:

First, suppose that $\iA v_{\frac{p-1}{2}} = c v_{\frac{p-1}{2}}$ for some $c 
\in \k.$ Since $\iA^2 = 1,$ we see that $c = \pm 1$.
Now let $$W = \span_{\k} \{ v_0, v_1, \ldots, v_{p-1}, w_{-\frac{k}{2}}, 
w_{-\frac{k}{2} + 1}, \ldots, w_{\frac{k}{2} - 1} \}.$$
From the above information, we can deduce that
\eqref{v14first}--\eqref{v14last} are satisfied.
Therefore, $V_{1,4}^{c} \subset
V.$ Since $V$ is irreducible, we conclude that $V = V_{1,4}^{c}$.
Since $\iA$ acts as $c \Id$ on $V_{1,4}^c \left[ \frac{p-1}{2} \right]$, we see 
that $V_{1,4}^{c}$ and $V_{1,4}^{c'}$ are isomorphic if and only if $c = c'$.
 
Second, suppose that $\iA v_{\frac{p-1}{2}}$ is not a scalar multiple of 
$v_{\frac{p-1}{2}}.$ In this case, we may write 
$\iA v_{\frac{p-1}{2}} = u_{\frac{p-1}{2}},$ and then
\begin{eqnarray*}
u_{\frac{p-1}{2}+j} &=& - \iB (\iA \iB)^{j-1} u_{\frac{p-1}{2}}, \quad j = 1, 
\ldots, \frac{p-1-k}{2} \\ 
u_{\frac{p-1}{2}-j} &=& (\iA \iB)^j u_{\frac{p-1}{2}} \quad j = 1,
\ldots, \frac{p-1-k}{2}.
\end{eqnarray*}
Since $\{ v_{\frac{p-1}{2}}, u_{\frac{p-1}{2}} \}$ is linearly independent,
it follows that $\{ v_j, u_j \}$ is linearly independent for all
$j = \frac{k}{2}, \frac{k}{2} + 1, \ldots, -\frac{k}{2} -1.$
Now let $$W = \span_{\k} \{v_0, v_1, \ldots, v_{p-1}, w_{-\frac{k}{2}}, 
w_{-\frac{k}{2}+1}, \ldots, w_{\frac{k}{2} - 1}, u_{\frac{k}{2}}, 
u_{\frac{k}{2}+1}, \ldots, u_{-\frac{k}{2}-1} \}.$$ 
Now, we know from lemma \ref{1lem8} that $\dm \vgen[-\frac{k}{2}] = 2,$ 
and this
forces $\iB u_{\frac{k}{2}} = c v_{-\frac{k}{2}}$ for some $c \in \k$.
From the above information, we can deduce that
\eqref{v15first}--\eqref{v15last} are satisfied.
Therefore, $V_{1,5}^{c} \subset
V.$ Since $V$ is irreducible, we conclude that $V = V_{1,5}^{c}$.
Since $\iB (\iA \iB)^{p-k}$ acts as $c \Id$ on $V_{1,5}^c \left[ -\frac{k}{2} 
\right]$, 
we see that $V_{1,5}^{c}$ and $V_{1,5}^{c'}$ are isomorphic if and only if $c = 
c'$.
\end{proof}

\begin{lemma}
Suppose $k = 0.$
Let $V$ be an irreducible representation of $\Hh_1$ on which the central
element $(\y^p - \y)^2$ acts as $0$. Then each eigenspace 
$V[c], c \in \pf$ has dimension at most $2$.
\label{1lem77}
\end{lemma}
                                                                                         
\begin{proof}
Let $\mu$ be an eigenvalue of $\y$, and let $v \in V[\mu].$ Note that $0
= (\y^p - \y)^2 v = (\mu^p - \mu)^2 v,$ so we have $\mu^p - \mu = 0.$
Hence, $\mu \in \pf.$
By corollary \ref{1cor4}, we have the following homomorphisms:
\begin{equation}
\vgen[0] \stackrel{iA}{\to} \vgen[-1]
\stackrel{\iB}{\to} \vgen[1]
\stackrel{\iA}{\to} \cdots
\stackrel{\iA}{\to} \vgen \left[ \frac{p+1}{2} \right]
\stackrel{\iB}{\to} \vgen \left[ \frac{p-1}{2} \right].
\label{eqn77}
\end{equation}
Since $k = 0$, corollary \ref{1cor4} implies that all of the homomorphisms in
\eqref{eqn77} are isomorphisms, and the eigenvalues of $\y$ are precisely
the elements of $\pf.$
Thus, $\dm V \geq  p \dm V [0]$.
Now the dimension of the algebra $\Hh_1 / (\X^p + \X^{-p} = a,
(\y^p - \y)^2 = 0)$ acting irreducibly on $V$ is at most $8 p^2$ (see the
proof of corollary \ref{1cor2}). Hence, $p^2 (\dm V[0])^2 \leq (\dm V)^2
\leq 8 p^2,$ which implies that $\dm V[0] \leq 2.$
\end{proof}
                                                                                         
\begin{proof}[Proof of Proposition \ref{prop5}]
As in the proof of Proposition \ref{prop3}, we see that $V_{1,1}^{\mu,d}$ is an
irreducible representation of $\Hh_1$ whenever $\mu$ and $d$ satisfy the conditions
in the statement of Proposition \ref{prop5}. Similarly, for
$c, \theta = \pm 1,$ $V_{1,6}^{c,\theta}$ is an irreducible
representation of $\Hh_1$ and for $c = \pm 1, a \in \k$, $V_{1,7}^{c,a}$ is an 
irreducible representation of $\Hh_1$.
                                                                                         
Now, let $V$ be an irreducible representation of $\Hh_1$,
and suppose that that $(\y^p - \y)^2$ acts on $V$ as $b$.
If $b \neq 0$, then $k$ is not a root of 
$f(y) = (y^p - y)^2 - b$, so the argument in the proof of 
Proposition \ref{prop3} applies for $V_{1,1}^{\mu, d}$. We will now assume 
that $b = 0$.
First we note that $\s \y = - \y \s$, which means that $\s V[0] \subset V[0].$
Now let $v_0 \in V[0]$ be an eigenvector of $\s$. Since $\s^2 = 1,$ we have
$\s v_0 = c v_0,$ where $c = \pm 1.$
Let
\begin{eqnarray*}
v_{-j} &=& \iA (\iB \iA)^{j-1} v_0, \quad j = 1, 2, \ldots, \frac{p-1}{2}; \\
v_j &=& (\iB \iA)^j v_0, \quad j = 1, 2, \ldots, \frac{p-1}{2}.
\end{eqnarray*}
 
Now there are two cases:

First, suppose that $\iA v_{\frac{p-1}{2}} = \theta 
v_{\frac{p-1}{2}}$ for some $\theta \in \k.$ 
Since $\iA^2 = 1,$ we see that $\theta = \pm 1$.
Now let $$W = \span_{\k} \{ v_0, v_1, \ldots, v_{p-1} \}.$$
From the above information, we can deduce that 
\eqref{v16first}--\eqref{v16last} are satisfied.
Therefore, $V_{1,6}^{c,\theta} \subset
V.$ Since $V$ is irreducible, we conclude that $V = V_{1,6}^{c,\theta}$.
Since $\s$ acts on $V[0]$ as $c \Id$ and $\iA$ acts on 
$V \left[\frac{p-1}{2} \right]$ as 
$\theta \Id$, we see that $V_{1,6}^{c,\theta}$ and $V_{1,6}^{c',\theta'}$ 
are isomorphic if and only if $c = c'$ and $\theta = \theta'$.

Second, suppose that $\iA v_{\frac{p-1}{2}}$ is not a scalar multiple of
$v_{\frac{p-1}{2}}.$ In this case, we may write
$\iA v_{\frac{p-1}{2}} = u_{\frac{p-1}{2}},$ and then
\begin{eqnarray*}
u_{\frac{p-1}{2}+j} &=& \iB (\iA \iB)^{j-1} u_{\frac{p-1}{2}}, \quad j = 1,
\ldots, \frac{p-3}{2} \\
u_{\frac{p-1}{2}-j} &=& (\iA \iB)^j u_{\frac{p-1}{2}} \quad j = 1,
\ldots, \frac{p-1}{2}.
\end{eqnarray*}
Since $\{ v_{\frac{p-1}{2}}, u_{\frac{p-1}{2}} \}$ is linearly independent,
it follows that $\{ v_j, u_j \}$ is linearly independent for all $j$. 
Hence, by lemma \ref{1lem77}, $V[0] = \span_{\k} \{ v_0, u_0 \},$ 
and this forces $\s u_0 = a v_0 + r u_0$ for some $a,r \in \k$.
Now $$u_0 = \s^2 u_0 = a \s v_0 + r \s u_0 = a c v_0 + a r v_0 + r^2 u_0.$$
Thus $r = \pm 1.$ Now if $c = r,$ then $a = 0$. But then $v_0 - u_0$ is an 
eigenvector of $\s$ and $\iA$ acts on $v_{\frac{p-1}{2}} - u_{\frac{p-1}{2}}$ 
as a scalar, and so the first case shows us that $V$ has a 
$p$-dimensional subrepresentation, contradicting $V$'s irreducibility. So we
may assume that $c = -r$; $a$ is then arbitrary. Now we can see that the 
eigenvalues of $\s$ acting on $V[0]$ are $\pm c$; that is, $\pm 1$. This 
means that we can assume $c=1$. From all this information, we can deduce 
that \eqref{v17first}--\eqref{v17last} are satisfied. Therefore, 
$V_{1,7}^a \subset V.$ Since $V$ is irreducible, we conclude that $V = 
V_{1,7}^a$.
Finally, $a$ is the coefficient of $v$ in
$\s \iA (\iB \iA)^{p-1} v$; thus, $V_{1,7}^a$ and $V_{1,7}^{a'}$ are 
isomorphic if and only if $a = a'$. 
\end{proof}

\bibliographystyle{alpha}
\bibliography{prob15}
\end{document}